# Generalized Convex Functions on Fractal Sets and Two Related Inequalities

**Huixia Mo**[1], **Xin Sui**[2] **and Dongyan Yu**[3]

[1,2,3]School of Science, Beijing University of Posts and Telecommunications, Beijing,100876, China,

Correspondence should be addressed to Hui-xia Mo; huixmo@bupt.edu.cn

In the paper, we introduce the generalized convex function on fractal sets $R^\alpha (0 < \alpha \leq 1)$ of real line numbers and study the properties of the generalized convex function. Based on these properties, we establish the generalized Jensen's inequality and generalized Hermite-Hadamard's inequality. Furthermore, some applications are given.

## 1. Introduction

Let $f: I \subseteq R \to R$. For any $x_1, x_2 \in I$ and $\lambda \in [0,1]$, if the following inequality

$$f(\lambda x_1 + (1-\lambda)x_2) \leq \lambda f(x_1) + (1-\lambda) f(x_2)$$

holds, then $f$ is called a convex function on $I$.

The convexity of functions play a significant role in many fields, for example in biological system, economy, optimization and so on [1-2]. And many important inequalities are established for the class of convex functions. For example, the Jensen's inequality and Hermite-Hadamard's inequality are the best known results in the literature, which can be stated as follows.

**Jensen's inequality** [3]: Assume that $f$ is a convex function on $[a,b]$. Then for any $x_i \in [a,b]$ and $\lambda_i \in [0,1]$ $(i = 1, 2, ..., n)$ with $\sum_{i=1}^{n} \lambda_i = 1$, we have

$$f\left(\sum_{i=1}^{n} \lambda_i x_i\right) \leq \sum_{i=1}^{n} \lambda_i f(x_i).$$

**Hermite-Hadamard's inequality** [4]: Let $f$ be a convex function on $[a,b]$ with $a < b$. If $f$ is integral on $[a,b]$, then

$$f\left(\frac{a+b}{2}\right) \leq \frac{1}{b-a} \int_a^b f(x)dx \leq \frac{f(a)+f(b)}{2}.$$

In recent years, the fractal has received significantly remarkable attention from scientists and engineers. In the sense of Mandelbrot, a fractal set is the one whose Hausdorff dimension strictly exceeds the topological dimension [5-9]. Many researchers studied the properties of functions on fractal space and constructed many kinds of fractional calculus by using different approaches (see [10-14]). Particularly, in [13], Yang stated the analysis of local fractional functions on fractal space systematically, which includes local fractional calculus, the monotonicity of function and so on.

Inspired by these investigations, we will introduce the generalized convex function on fractal sets and establish the generalized Jensen's inequality and generalized Hermite-Hadamard's inequality related to generalized convex function. We shall focus our attention on the convexity since a function $f$ is concave if and only if $-f$ is convex. So, every result for the convex function can be easily re-stated in terms of concave functions.

The article is organized as follows: In Section 2, we state the operations with real line number on fractal sets and give the definitions of the local fractional derivatives and local fractional integral. In Section 3, we introduce the definition of the generalized convex function on fractal sets and study the properties of the generalized convex functions. In Section 4, we establish the generalized Jensen's inequality and generalized Hermite-Hadamard's inequality on fractal sets. In Section 5, some applications are given on fractal sets by means of the generalized Jensen's inequality.

## 2. Preliminaries

Recall the set $R^\alpha$ of real line numbers and use the Gao-Yang-Kang's idea to describe the definitions of the local fractional derivative and local fractional integral.

Recently, the theory of Yang's fractional sets [13] was introduced as follows.

For $0 < \alpha \leq 1$, we have the following $\alpha$-type set of element sets:

$Z^\alpha$: The $\alpha$-type set of the integer are defined as the set $\{0^\alpha, \pm 1^\alpha, \pm 2^\alpha, \cdots, \pm n^\alpha, \cdots\}$.

$Q^\alpha$: The $\alpha$-type set of the rational numbers are defined as the set $\{m^\alpha = (p/q)^\alpha : p \in Z, q \neq 0\}$.

$J^\alpha$: The $\alpha$-type set of the irrational numbers are defined as the set $\{m^\alpha \neq (p/q)^\alpha : p \in Z, q \neq 0\}$.

$R^\alpha$: The $\alpha$-type set of the real line numbers are defined as the set $R^\alpha = Q^\alpha \bigcup J^\alpha$.

If $a^\alpha, b^\alpha$ and $c^\alpha$ belong to the set $R^\alpha$ of real line numbers, then

(1) $a^\alpha + b^\alpha$ and $a^\alpha b^\alpha$ belong to the set $R^\alpha$;

(2) $a^\alpha + b^\alpha = b^\alpha + a^\alpha = (a+b)^\alpha = (b+a)^\alpha$;

(3) $a^\alpha + (b^\alpha + c^\alpha) = (a+b)^\alpha + c^\alpha$;

(4) $a^\alpha b^\alpha = b^\alpha a^\alpha = (ab)^\alpha = (ba)^\alpha$;

(5) $a^\alpha (b^\alpha c^\alpha) = (a^\alpha b^\alpha) c^\alpha$;

(6) $a^\alpha (b^\alpha + c^\alpha) = a^\alpha b^\alpha + a^\alpha c^\alpha$;

(7) $a^\alpha + 0^\alpha = 0^\alpha + a^\alpha = a^\alpha$ and $a^\alpha 1^\alpha = 1^\alpha a^\alpha = a^\alpha$.

Let us now state some definitions about the local fractional calculus on $R^\alpha$.

**Definition 2.1** [13] A non-differentiable function $f : R \to R^\alpha, x \to f(x)$ is called to be local fractional continuous at $x_0$, if for any $\varepsilon > 0$, there exists $\delta > 0$, such that

$$|f(x) - f(x_0)| < \varepsilon^\alpha$$

holds for $|x - x_0| < \delta$, where $\varepsilon, \delta \in R$. If $f(x)$ is local fractional continuous on the interval $(a, b)$, we denote $f(x) \in C_\alpha(a, b)$.

**Definition 2.2** [13] The local fractional derivative of $f(x)$ of order $\alpha$ at $x = x_0$ is defined

by

$$f^{(\alpha)}(x_0) = \left.\frac{d^\alpha f(x)}{dx^\alpha}\right|_{x=x_0} = \lim_{x \to x_0} \frac{\Delta^\alpha(f(x)-f(x_0))}{(x-x_0)^\alpha},$$

where $\Delta^\alpha(f(x)-f(x_0)) = \Gamma(1+\alpha)(f(x)-f(x_0))$.

If there exists $f^{((k+1)\alpha)}(x) = \overbrace{D_x^\alpha \cdots D_x^\alpha}^{k+1 \text{ times}} f(x)$ for any $x \in I \subseteq R$, then we denote $f \in D_{(k+1)\alpha}(I)$, where $k = 0, 1, 2, \ldots$.

**Definition 2.3** [13] The local fractional integral of the function $f(x)$ of order $\alpha$ is defined by

$$_aI_b^{(\alpha)} f(x)$$
$$= \frac{1}{\Gamma(1+a)} \int_a^b f(t)(dt)^\alpha$$
$$= \frac{1}{\Gamma(1+a)} \lim_{\Delta t \to 0} \sum_{j=0}^{N-1} f(t_j)(\Delta t_j)^\alpha,$$

with $\Delta t_j = t_{j+1} - t_j$ and $\Delta t = \max\{\Delta t_j \mid j = 1, 2, \ldots, N-1\}$, where $[t_j, t_{j+1}]$, $j = 0, \ldots, N-1$ and $t_0 = a < t_1 < \cdots < t_i < \cdots < t_{N-1} < t_N = b$ is a partition of the interval $[a,b]$.

Here, it follows that $_aI_a^{(\alpha)} f(x) = 0$ if $a = b$ and $_aI_b^{(\alpha)} f(x) = -_bI_a^{(\alpha)} f(x)$ if $a < b$. If for any $x \in [a,b]$, there exists $_aI_x^{(\alpha)} f(x)$, then it is denoted by $f(x) \in I_x^{(\alpha)}[a,b]$.

**Lemma 2.1** [13] (Generalized local fractional Taylor theorem)
Suppose that $f^{(k+1)\alpha}(x) \in C_\alpha(I)$, for interval $I \subseteq R$, $k = 0, 1 \ldots n$, $0 < \alpha \le 1$. And let $x_0 \in [a,b]$. Then for any $x \in I$, there exists at least one point $\xi$, which lies between the points $x$ and $x_0$, such that

$$f(x) = \sum_{k=0}^{n} \frac{f^{(k\alpha)}(x_0)}{\Gamma(1+ka)}(x-x_0)^{k\alpha} + \frac{f^{((n+1)\alpha)}(\xi)}{\Gamma(1+(n+1)a)}(x-x_0)^{(n+1)\alpha}.$$

**Remark 2.1** When $I \subseteq R$ is an open interval $(a,b)$, Yang [13] has given the proof for the generalized local fractional Taylor theorem. In fact, using the generalized local fractional Lagrange's theorem and following the proof of the class Taylor theorem, we can show that for any interval $I \subseteq R$, the formula is also true.

## 3. Generalized convex functions

From an analytical point of view, we have the following definition.

**Definition 3.1** Let $f : I \subseteq R \to R^\alpha$. For any $x_1, x_2 \in I$ and $\lambda \in [0,1]$, if the following inequality

$$f(\lambda x_1 + (1-\lambda)x_2) \le \lambda^\alpha f(x_1) + (1-\lambda)^\alpha f(x_2)$$

holds, then $f$ is called a generalized convex function on $I$.

**Definition 3.2** Let $f : I \to R^\alpha$. For any $x_1 \ne x_2 \in I$ and $\lambda \in [0,1]$, if the following inequality

$$f(\lambda x_1 + (1-\lambda)x_2) < \lambda^\alpha f(x_1) + (1-\lambda)^\alpha f(x_2)$$

holds, then $f$ is called a generalized strictly convex function on $I \subseteq R$.

It follows immediately, from the given definitions, that a generalized strictly convex function is also generalized convex. But, the converse is not true. And if these two inequalities is reversed, then $f$ is called a generalized concave function or generalized strictly concave function, respectively.

Here are two basic examples of generalized strictly convex functions:

(1) $f(x) = x^{\alpha p}$, $x \geq 0$, $p > 1$;

(2) $f(x) = E_\alpha(x^\alpha)$, $x \in R$, where $E_\alpha(x^\alpha) = \sum_{k=0}^{\infty} \frac{x^{\alpha k}}{\Gamma(1+k\alpha)}$ is the Mittag-Leffer function.

Note that the linear function $f(x) = a^\alpha x^\alpha + b^\alpha$, $x \in R$ is generalized convex and also generalized concave.

We shall focus our attention on the convexity since a function $f$ is concave if and only if $-f$ is convex. So, every result for the convex function can be easily re-stated in terms of concave functions.

In the following, we will study the properties of the generalized convex functions.

**Theorem 3.1** Let $f: I \to R^\alpha$. Then $f$ is a generalized convex function if and only if the inequality

$$\frac{f(x_1) - f(x_2)}{(x_1 - x_2)^\alpha} \leq \frac{f(x_3) - f(x_2)}{(x_3 - x_2)^\alpha}$$

holds, for any $x_1, x_2, x_3 \in I$ with $x_1 < x_2 < x_3$.

**Proof.** In fact, take $\lambda = \frac{x_3 - x_2}{x_3 - x_1}$, then $x_2 = \lambda x_1 + (1-\lambda) x_3$. And by the generalized convexity of $f$, we get

$$f(x_2) = f(\lambda x_1 + (1-\lambda) x_3) \leq \lambda^\alpha f(x_1) + (1-\lambda)^\alpha f(x_3) = \left(\frac{x_3 - x_2}{x_3 - x_1}\right)^\alpha f(x_1) + \left(\frac{x_2 - x_1}{x_3 - x_1}\right)^\alpha f(x_3).$$

From the above formula, it is easy to see that

$$\frac{f(x_1) - f(x_2)}{(x_1 - x_2)^\alpha} \leq \frac{f(x_3) - f(x_2)}{(x_3 - x_2)^\alpha}.$$

Reversely, for any two points $x_1, x_3 (x_1 < x_3)$ on $I \subseteq R$, we take $x_2 = \lambda x_1 + (1-\lambda) x_3$ for $\lambda \in (0,1)$. Then $x_1 < x_2 < x_3$ and $\lambda = \frac{x_3 - x_2}{x_3 - x_1}$. Using the above inverse process, we have

$$f(\lambda x_1 + (1-\lambda) x_3) \leq \lambda^\alpha f(x_1) + (1-\lambda)^\alpha f(x_3).$$

So, $f$ is a convex function on $I \subseteq R$.

In the same way, it can be shown that $f$ is a generalized convex function on $I \subseteq R$ if and only if

$$\frac{f(x_2) - f(x_1)}{(x_2 - x_1)^\alpha} \leq \frac{f(x_3) - f(x_1)}{(x_3 - x_1)^\alpha} \leq \frac{f(x_3) - f(x_2)}{(x_3 - x_2)^\alpha},$$

for any $x_1, x_2, x_3 \in I$ with $x_1 < x_2 < x_3$.

**Theorem 3.2** Let $f \in D_\alpha(I)$, then the following conditions are equivalent.

(1) $f$ is a generalized convex function on $I$,

(2) $f^{(\alpha)}$ is an increasing function on $I$,

(3) for any $x_1, x_2 \in I$,

$$f(x_2) \geq f(x_1) + \frac{f^{(\alpha)}(x_1)}{\Gamma(1+\alpha)}(x_2 - x_1)^\alpha.$$

**Proof.** $(1 \to 2)$ Let $x_1, x_2 \in I$ with $x_1 < x_2$. And take $h > 0$ which is small enough such that $x_1 - h, x_2 + h \in I$. Since $x_1 - h < x_1 < x_2 < x_2 + h$, then using Theorem 3.1 we have

$$\Gamma(1+a)\frac{f(x_1) - f(x_1 - h)}{h^\alpha} \leq \Gamma(1+a)\frac{f(x_2) - f(x_1)}{(x_2 - x_1)^\alpha} \leq \Gamma(1+a)\frac{f(x_2 + h) - f(x_2)}{h^\alpha}.$$

Since $f \in D_\alpha(I)$, then let $h \to 0^+$, it follows that

$$f^{(\alpha)}(x_1) \leq \Gamma(1+a)\frac{f(x_2) - f(x_1)}{(x_2 - x_1)^\alpha} \leq f^{(\alpha)}(x_2).$$

So, $f^{(\alpha)}$ is increasing in $I$.

$(2 \to 3)$ Take $x_1, x_2 \in I$. Without loss of generality, we can assume that $x_1 < x_2$. Since $f^{(\alpha)}$ is increasing in the interval $I$, then applying the generalized local fractional Taylor theorem, we have

$$f(x_2) - f(x_1) = \frac{f^{(\alpha)}(\xi)}{\Gamma(1+a)}(x_2 - x_1)^\alpha \geq \frac{f^{(\alpha)}(x_1)}{\Gamma(1+a)}(x_2 - x_1)^\alpha,$$

where $\xi \in (x_1, x_2)$. That is to say

$$f(x_2) \geq f(x_1) + \frac{f^{(\alpha)}(x_1)}{\Gamma(1+a)}(x_2 - x_1)^\alpha.$$

$(3 \to 1)$ For any $x_1, x_2 \in I$, we let $x_3 = \lambda x_1 + (1-\lambda)x_2$, where $0 < \lambda < 1$. It is easy to see that $x_1 - x_3 = (1-\lambda)(x_1 - x_2)$ and $x_2 - x_3 = \lambda(x_2 - x_1)$. Then from the third condition, we have

$$f(x_1) \geq f(x_3) + \frac{f^{(\alpha)}(x_3)}{\Gamma(1+a)}(x_1 - x_3)^\alpha = f(x_3) + (1-\lambda)^\alpha \frac{f^{(\alpha)}(x_3)}{\Gamma(1+a)}(x_1 - x_2)^\alpha,$$

and

$$f(x_2) \geq f(x_3) + \frac{f^{(\alpha)}(x_3)}{\Gamma(1+a)}(x_2 - x_3)^\alpha = f(x_3) + \lambda^\alpha \frac{f^{(\alpha)}(x_3)}{\Gamma(1+a)}(x_2 - x_1)^\alpha.$$

At the above two formulas, multiply $\lambda^\alpha$ and $(1-\lambda)^\alpha$, respectively, then we obtain

$$\lambda^\alpha f(x_1) + (1-\lambda)^\alpha f(x_2) \geq f(x_3) = f(\lambda x_1 + (1-\lambda)x_2).$$

So, $f$ is a generalized convex function on $I$.

**Corollary 3.1** Let $f \in D_{2\alpha}(a,b)$. Then $f$ is a generalized convex function (or a generalized concave function) if and only if

$$f^{(2\alpha)}(x) \geq 0 (\text{or } f^{(2\alpha)}(x) \leq 0),$$

for any $x \in (a,b)$.

## 4. Some inequalities

**Theorem 4.1 (Generalized Jensen's inequality)** Assume that $f$ is a generalized convex function on $[a,b]$. Then for any $x_i \in [a,b]$ and $\lambda_i \in [0,1]$ $(i=1, 2,..., n)$ with $\sum_{i=1}^{n} \lambda_i = 1$, we have

$$f\left(\sum_{i=1}^{n} \lambda_i x_i\right) \leq \sum_{i=1}^{n} \lambda_i^{\alpha} f(x_i).$$

**Proof.** When $n=2$, the inequality is obviously true. Assume that for $n=k$ the inequality is also true. Then for any $x_1, x_2, ..., x_k \in [a,b]$ and $\gamma_i > 0$, $i=1,2...,k$ with $\sum_{i=1}^{k} \gamma_i = 1$, we have

$$f\left(\sum_{i=1}^{k} \gamma_i x_i\right) \leq \sum_{i=1}^{k} \gamma_i^{\alpha} f(x_i).$$

If $x_1, x_2, ..., x_k, x_{(k+1)} \in [a,b]$ and $\lambda_i > 0$ for $i=1,2,...,k+1$ with $\sum_{i=1}^{k+1} \lambda_i = 1$, then one sets up $\gamma_i = \dfrac{\lambda_i}{1-\lambda_{k+1}}$, $i=1,2,...,k$. It is easy to see $\sum_{i=1}^{k} \gamma_i = 1$.

Thus,
$$f(\lambda_1 x_1 + \lambda_2 x_2 + ... + \lambda_k x_k + \lambda_{k+1} x_{k+1})$$
$$= f\left((1-\lambda_{k+1})\frac{\lambda_1 x_1 + \lambda_2 x_2 + ... + \lambda_k x_k}{1-\lambda_{k+1}} + \lambda_{k+1} x_{k+1}\right)$$
$$\leq (1-\lambda_{k+1})^{\alpha} f(\gamma_1 x_1 + \gamma_2 x_2 + ... + \gamma_k x_k) + \lambda_{k+1}^{\alpha} f(x_{k+1})$$
$$\leq (1-\lambda_{k+1})^{\alpha} [\gamma_1^{\alpha} f(x_1) + \gamma_2^{\alpha} f(x_2) + ... + \gamma_k^{\alpha} f(x_k)] + \lambda_{k+1}^{\alpha} f(x_{k+1})$$
$$= (1-\lambda_{k+1})^{\alpha} \left[\left(\frac{\lambda_1}{1-\lambda_{k+1}}\right)^{\alpha} f(x_1) + \left(\frac{\lambda_2}{1-\lambda_{k+1}}\right)^{\alpha} f(x_2) + ... + \left(\frac{\lambda_k}{1-\lambda_{k+1}}\right)^{\alpha} f(x_k)\right] + \lambda_{k+1}^{\alpha} f(x_{k+1})$$
$$= \sum_{i=1}^{k} \lambda_i^{\alpha} f(x_i).$$

So, the mathematical induction gives the proof of Theorem 4.1.

**Corollary 4.1** Let $f \in D_{2\alpha}[a,b]$ and $f^{(2\alpha)}(x) \geq 0$ for any $x \in [a,b]$. Then for any $x_i \in [a,b]$ and $\lambda_i \in [0,1]$ $(i=1,2,...n)$ with $\sum_{i=1}^{n} \lambda_i = 1$, we have

$$f\left(\sum_{i=1}^{n} \lambda_i x_i\right) \leq \sum_{i=1}^{n} \lambda_i^{\alpha} f(x_i).$$

Using the generalized Jensen's inequality and the convexity of functions, we can also get some integral inequalities.

In [13], Yang established the generalized Cauchy-Schwatz's inequality by the estimate $a^{\frac{\alpha}{p}} b^{\frac{\alpha}{q}} \leq \dfrac{a^{\alpha}}{p} \dfrac{b^{\alpha}}{q}$, where $a^{\alpha}, b^{\alpha} > 0$, $p, q \geq 1$ and $\dfrac{1}{p} + \dfrac{1}{q} = 1$.

Now, via the generalized Jensen's inequality, we will give another proof for the generalized Cauchy-Schwarz's inequality.

**Corollary 4.2 (Generalized Cauchy-Schwarz's inequality)** Let $|a_k|>0$, $|b_k|>0$, $k = 1, 2,...,n$. Then we have

$$\sum_{k=1}^{n}|a_k|^{\alpha}|b_k|^{\alpha} \leq \left(\sum_{k=1}^{n}|a_k|^{2\alpha}\right)^{\frac{1}{2}}\left(\sum_{k=1}^{n}|b_k|^{2\alpha}\right)^{\frac{1}{2}}.$$

**Proof.** Take $f(x) = x^{2\alpha}$. It is easy to see that $f^{(2\alpha)}(x) \geq 0$ for any $x \in (a,b)$.

Take

$$\lambda_k = \frac{|b_k|^2}{\sum_{k=1}^{n}|b_k|^2}, \quad x_k = \frac{|a_k|}{|b_k|}. \text{ Then } 0 \leq \lambda_k \leq 1 \ (k=1,2,...,n) \text{ with } \sum_{k=1}^{n}\lambda_k = 1.$$

Thus, by Jensen's inequality $f\left(\sum_{k=1}^{n}\lambda_k x_k\right) \leq \sum_{k=1}^{n}\lambda_k^{(\alpha)} f(x_k)$, we have

$$\left[\sum_{k=1}^{n}\frac{|b_k|^2}{\sum_{k=1}^{n}|b_k|^2}\frac{|a_k|}{|b_k|}\right]^{2\alpha} \leq \sum_{k=1}^{n}\left[\frac{|b_k|^2}{\sum_{k=1}^{n}|b_k|^2}\right]^{\alpha}\left[\frac{|a_k|}{|b_k|}\right]^{2\alpha}.$$

The above formula can be reduced to

$$\left[\sum_{k=1}^{n}\frac{|b_k||a_k|}{\sum_{k=1}^{n}|b_k|^2}\right]^{2\alpha} \leq \sum_{k=1}^{n}\left[\frac{|a_k|^{2\alpha}}{(\sum_{k=1}^{n}|b_k|^2)^{\alpha}}\right],$$

which implies that

$$\left[\sum_{k=1}^{n}|b_k||a_k|\right]^{2\alpha} \leq \sum_{k=1}^{n}|a_k|^{2\alpha}\sum_{k=1}^{n}|b_k|^{2\alpha}.$$

Thus we have

$$\sum_{k=1}^{n}|a_k|^{\alpha}|b_k|^{\alpha} \leq \left(\sum_{k=1}^{n}|a_k|^{2\alpha}\right)^{\frac{1}{2}}\left(\sum_{k=1}^{n}|b_k|^{2\alpha}\right)^{\frac{1}{2}}.$$

**Theorem 4.2 (Generalized Hermite-Hadmard's inequality)** Let $f(x) \in I_x^{(\alpha)}[a,b]$ be a generalized convex function on $[a,b]$ with $a < b$. Then

$$f\left(\frac{a+b}{2}\right) \leq \frac{\Gamma(1+\alpha)}{(b-a)^{\alpha}} {}_aI_b^{(\alpha)} f(x) \leq \frac{f(a)+f(b)}{2^{\alpha}}.$$

**Proof.** Let $x = a+b-y$. Then

$$\int_a^{\frac{a+b}{2}} f(x)(dx)^\alpha = \int_{\frac{a+b}{2}}^b f(a+b-y)(dy)^\alpha.$$

Furthermore, when $x \in \left[\frac{a+b}{2}, b\right]$, $a+b-x \in \left[a, \frac{a+b}{2}\right]$. And by the convexity of $f$, we have

$$f(a+b-x) + f(x) \geq 2^\alpha f\left(\frac{a+b}{2}\right).$$

Thus

$$\int_a^b f(x)(dx)^\alpha$$
$$= \int_a^{\frac{a+b}{2}} f(x)(dx)^\alpha + \int_{\frac{a+b}{2}}^b f(x)(dx)^\alpha$$
$$= \int_{\frac{a+b}{2}}^b [f(a+b-x) + f(x)](dx)^\alpha$$
$$\geq \int_{\frac{a+b}{2}}^b 2^\alpha f\left(\frac{a+b}{2}\right)(dx)^\alpha$$
$$= (b-a)^\alpha f\left(\frac{a+b}{2}\right). \tag{4.1}$$

For another part, we first note that if $f$ is a generalized convex function, then, for $t \in [0,1]$, it yields

$$f(ta + (1-t)b) \leq t^\alpha f(a) + (1-t)^\alpha f(b),$$

and

$$f((1-t)a + tb) \leq (1-t)^\alpha f(a) + t^\alpha f(b).$$

By adding these inequalities we have

$$f(ta+(1-t)b) + f((1-t)a+tb) \leq t^\alpha f(a) + (1-t)^\alpha f(b) + (1-t)^\alpha f(a) + t^\alpha f(b) = f(a) + f(b).$$

Then, integrating the resulting inequality with respect to $t$ over $[0,1]$, we obtain

$$\frac{1}{\Gamma(1+\alpha)} \int_0^1 [f(ta+(1-t)b) + f((1-t)a+tb)](dt)^\alpha$$
$$\leq \frac{1}{\Gamma(1+\alpha)} \int_0^1 (f(a)+f(b))(dt)^\alpha.$$

It is easy to see that

$$\frac{1}{\Gamma(1+\alpha)} \int_0^1 [f(ta+(1-t)b) + f((1-t)a+tb)](dt)^\alpha = \frac{2^\alpha}{(b-a)^\alpha} {}_a I_b^{(\alpha)} f(x),$$

and

$$\frac{1}{\Gamma(1+\alpha)} \int_0^1 (f(a)+f(b))(dt)^\alpha = \frac{f(a)+f(b)}{\Gamma(1+\alpha)}.$$

So,

$$\frac{\Gamma(1+\alpha)}{(b-a)^\alpha} {}_aI_b^{(\alpha)} f(x) \leq \frac{f(a)+f(b)}{2^\alpha}. \tag{4.2}$$

Combing the inequalities (4.1) and (4.2), we have

$$f\left(\frac{a+b}{2}\right) \leq \frac{\Gamma(1+\alpha)}{(b-a)^\alpha} {}_aI_b^{(\alpha)} f(x) \leq \frac{f(a)+f(b)}{2^\alpha}.$$

Note that, it will be reduced to the class Hermite- Hadmard inequality if $\alpha =1$.

## 5. Applications of generalized Jensen's inequality

Using the generalized Jensen's inequality, we can get some inequalities.

**Example 5.1** Let $a>0$, $b>0$ and $a^{3\alpha}+b^{3\alpha} \leq 2^\alpha$. Then $a+b \leq 2$.

**Proof.** Let $f(x)=x^{3\alpha}$, $x\in(0,+\infty)$. It is easy to see that $f$ is a generalized convex function. So,

$$f\left(\frac{a+b}{2}\right) \leq \frac{f(a)+f(b)}{2^\alpha}.$$

That is

$$\frac{(a+b)^{3\alpha}}{8^\alpha} \leq \frac{a^{3\alpha}+b^{3\alpha}}{2^\alpha} \leq 1^\alpha.$$

Thus, we conclude that $a+b \leq 2$.

**Example 5.2** Let $x, y \in R$. Then

$$E_\alpha\left(\left(\frac{x+y}{2}\right)^\alpha\right) \leq \frac{1}{2^\alpha}(E_\alpha(x^\alpha)+E_\alpha(y^\alpha)),$$

Where $E_\alpha(x^\alpha) = \sum_{k=0}^{\infty} \frac{x^{\alpha k}}{\Gamma(1+k\alpha)}$ is the Mittag-Leffer function.

**Proof.** Take $f(x)=E_\alpha(x^\alpha)$. It is easy to see $(E_\alpha(x^\alpha))^{(2\alpha)} = E_\alpha(x^\alpha) > 0$. So, the generalized Jensen's inequality gives

$$E_\alpha\left(\left(\frac{x+y}{2}\right)^\alpha\right) \leq \frac{1}{2^\alpha}(E_\alpha(x^\alpha)+E_\alpha(y^\alpha)).$$

**Example 5.3** (Power Mean Inequality) Let $a_1, a_2, \cdots a_n > 0$ and $0<s<t$ or $s<t<0$. Denote

$$S_r = \left(\frac{a_1^{\alpha r} + a_2^{\alpha r} + \cdots + a_n^{\alpha r}}{n^\alpha}\right)^{1/r}, \quad r\in R,$$

Then $S_s \leq S_t$. And $S_s = S_t$ if and only if $a_1 = a_2 = \cdots = a_n$.

**Proof. Case I:** $0<s<t$.

Take $f(x)=x^{(t/s)\alpha}$, $x>0$. Then

$$f^{(2\alpha)}(x) = \frac{\Gamma(1+\frac{t\alpha}{s})}{\Gamma(1+(\frac{t}{s}-1)\alpha)} x^{(t/s-2)\alpha} > 0.$$

By the generalized Jensen's inequality, we have

$$f\left(\frac{a_1^s + a_2^s + \cdots + a_n^s}{n}\right) \le \frac{f(a_1^s) + f(a_2^s) + \cdots + f(a_n^s)}{n^\alpha}.$$

That is

$$\left(\frac{a_1^s + a_2^s + \cdots + a_n^s}{n}\right)^{(t/s)\alpha} \le \frac{(a_1^s)^{(t/s)\alpha} + (a_2^s)^{(t/s)\alpha} + \cdots + (a_n^s)^{(t/s)\alpha}}{n^\alpha}.$$

From the above formula, it is easy to see

$$\left(\frac{a_1^{\alpha s} + a_2^{\alpha s} + \cdots + a_n^{\alpha s}}{n^\alpha}\right)^{1/s} \le \left(\frac{a_1^{\alpha t} + a_2^{\alpha t} + \cdots + a_n^{\alpha t}}{n^\alpha}\right)^{1/t}.$$

So, we have $S_s \le S_t$.

**Case II:** $s < t < 0$.

Let $b_i = 1/a_i$ and apply the case for $0 < -t < -s$, we can get the conclusion.

**Example 5.4** If $a, b, c > 0$ and $a + b + c = 1$, then find the minimum of

$$\left(a + \frac{1}{a}\right)^{10\alpha} + \left(b + \frac{1}{b}\right)^{10\alpha} + \left(c + \frac{1}{c}\right)^{10\alpha}.$$

**Solution.** Note that $0 < a, b, c < 1$. Let $f(x) = \left(x + \frac{1}{x}\right)^{10\alpha}$, $x \in (0,1)$. Then, via the formula

$$\frac{d^\alpha x^{k\alpha}}{dx^\alpha} = \frac{\Gamma(1+k\alpha)}{\Gamma(1+(k-1)\alpha)} x^{(k-1)\alpha},$$

we have

$$f^{(2\alpha)}(x) = \frac{\Gamma(1+10\alpha)}{\Gamma(1+8\alpha)}\left(x + \frac{1}{x}\right)^{8\alpha}\left(1 - \frac{1}{x^2}\right)^{2\alpha} + \frac{\Gamma(1+\alpha)\Gamma(1+10\alpha)}{\Gamma(1+9\alpha)}\left(x + \frac{1}{x}\right)^{9\alpha}\left(\frac{2}{x^3}\right)^\alpha > 0.$$

By the generalized Jensen's inequality,

$$\left(\frac{10}{3}\right)^{10\alpha} = f\left(\frac{a+b+c}{3}\right)$$

$$\le \frac{1}{3^\alpha}[f(a) + f(b) + f(c)]$$

$$= \frac{1}{3^\alpha}\left[\left(a + \frac{1}{a}\right)^{10\alpha} + \left(b + \frac{1}{b}\right)^{10\alpha} + \left(c + \frac{1}{c}\right)^{10\alpha}\right].$$

So, the minimum is $\dfrac{10^{10\alpha}}{3^{9\alpha}}$, when $a=b=c=1/3$.

**Example 5.5** If $a,b,c,d > 0$ and $c^{2\alpha}+d^{2\alpha}=(a^{2\alpha}+b^{2\alpha})^3$, then show that

$$\frac{a^{3\alpha}}{c^{\alpha}}+\frac{b^{3\alpha}}{d^{\alpha}} \geq 1^{\alpha}.$$

**Proof.** Let $x_1 = \left(\dfrac{a^3}{c}\right)^{1/2}$, $x_2 = \left(\dfrac{b^3}{d}\right)^{1/2}$, $y_1 = (ac)^{1/2}$, $y_2 = (bd)^{1/2}$. By the generalized Cauchy-Schwartz inequality, we have

$$\left(\frac{a^{3\alpha}}{c^{\alpha}}+\frac{b^{3\alpha}}{d^{\alpha}}\right)(a^{\alpha}c^{\alpha}+b^{\alpha}d^{\alpha})$$
$$= (x_1^{2\alpha}+x_2^{2\alpha})(y_1^{2\alpha}+y_2^{2\alpha})$$
$$\geq (x_1^{\alpha}y_1^{\alpha}+x_2^{\alpha}y_2^{\alpha})^2$$
$$= (a^{2\alpha}+b^{2\alpha})^2$$
$$= (a^{2\alpha}+b^{2\alpha})^{1/2}(c^{2\alpha}+d^{2\alpha})^{1/2}$$
$$\geq a^{\alpha}c^{\alpha}+b^{\alpha}d^{\alpha}.$$

Canceling $a^{\alpha}c^{\alpha}+b^{\alpha}d^{\alpha}$ on both sides, we get the desired result.

**5. Conclusion**

In the paper, we introduce the definition of generalized convex function on fractal sets. Based on the definition, we study the properties of the generalized convex functions and establish two important inequalities: the generalized Jensen's inequality and generalized Hermite-Hadamard's inequality. At last, we also give some applications for these inequalities on fractal sets.

**Acknowledgments**

The authors would like to express their gratitude to the reviewers for their very valuable comments. And, this work is supported by the National Natural Science Foundation of China (No. 11161042).